\documentclass [12pt] {article}
\usepackage{amsfonts}
\usepackage{amssymb}

\newcommand{\qed} {\hspace {0.1in} \rule {1.5mm} {3.5mm}}

\newtheorem{lemma}{Lemma}[section]

\newtheorem{theorem}{Theorem}

\newtheorem{proposition}{Proposition}[section]

\def\limn{\lim_{n\to\infty}}

\def\e{\epsilon}

\def\d{\delta}

\def\deg{\mbox{deg}\,}

\def\limo{\lim_{\omega}}

\def\Fo{F$\mbox{\o}$lner}

\def\proof{\smallskip\noindent{\it Proof.} }

\def\bZ{{\mathbb Z}}
\def\bR{{\mathbb R}}
\def\bN{{\mathbb N}}

\def\mi{^{(m)}}
\def\deg{\mbox{deg}\,}

\def\cG{\mbox{$\cal G$}}

\def\to{\rightarrow}

\def\uw{\underline{w}}

\def\wp{\widetilde{\Phi}}

\begin{document}
\title{Minimal topological actions do not determine the measurable
orbit equivalence class}
\author{\sc Tullio Ceccherini-Silberstein \footnote{Dipartimento di Ingegneria, Universit\`a del Sannio,
  C.so Garibaldi 108, 82100 Benevento, Italy. Email: {\tt tceccher@mat.uniroma1.it}.} and
G\'abor Elek
\footnote {The Alfred Renyi Mathematical Institute of
the Hungarian Academy of Sciences, P.O. Box 127, H-1364 Budapest, Hungary.
Email: {\tt elek@renyi.hu}, Supported by OTKA Grants T 049841 and T 037846, and by GNAMPA-INdAM.}}
\date{\today}
\maketitle \vskip 0.2in \noindent{\bf Abstract.} We construct
a minimal topological action $\wp$ of a non-amenable group
on a Cantor set $C$, which is non-uniquely ergodic and furthermore
there exist ergodic invariant measures $\mu_1$ and $\mu_2$
such that $(\wp,C,\mu_1)$ and $(\wp,C,\mu_2)$ are not orbit equivalent
measurable equivalence relations.
amenable
graphs 
\vskip 0.3in
\newpage
\section{Introduction}
\subsection{The result}
In \cite{OX} Oxtoby constructed a non-uniquely ergodic minimal
$\bZ$-action on a compact metric space
. It is well-known however that any two measure-preserving
transformations of countably infinite amenable groups are orbit equivalent.
Hence if $\mu_1$ and $\mu_2$ are ergodic invariant measures on a set $X$
under a $\bZ$-action $\wp$, then $(\wp,X,\mu_1)$ and $(\wp,X,\mu_2)$
are necessarily orbit equivalent. In this paper we consider the actions
of non-amenable groups. 

\noindent
Let $\Gamma$ be a countable non-amenable group, $X$ be a compact metric
space and \\ $\wp:\Gamma\to Homeo(X)$ be an injective homomorphism, in other words
$\wp$ is a faithful topological action of $\Gamma$ on $X$. Suppose that 
for any $x\in X$ the orbit $Orb(x)$ is dense, that is $\wp$ is minimal.
We consider such topological actions $\wp$, for which there exist
$\Gamma$-invariant measures on $X$ (for amenable groups such measures
always exist). The goal of our paper
 is to construct an action $\wp$ with invariant
ergodic measures $\mu_1$ and $\mu_2$ such that the measurable equivalence
relations $(\wp,X,\mu_1)$ and $(\wp,X,\mu_2)$ are not orbit equivalent.

\subsection{Gromov's graphs of dense holonomy}
Let $G(V,E)$ be an infinite connected graph with edges colored by the
symbols
 $A,B,C$ and $D$. Note that we consider only proper edge-colorings,
that is for each vertex $p\in V$ the colors of the edges
incident to $p$ are different. In particular, $G$ is
    of bounded valence 4 (each vertex has at most 4 neighbours). We call these graphs $4$-colored graphs.
Any such colorings define an action $\Phi$ of the group $\Gamma=\bZ_2*\bZ_2*\bZ_2*
\bZ_2$, the free product of four copies of the (cyclic) group of order two, on the
vertices of $G$. Indeed, for any word $\uw=w_nw_{n-1}\dots w_1$ in $\Gamma$ and $x \in V$, there
exists a unique path $(x_0,x_1,\dots,x_n)$ such that
\begin{itemize}
\item $x_0=x$
\item If $x_i$ has no incident edges colored by $w_i$, then $x_{i+1}=x_i$.
\item If the edge $(x_i,y)$ is colored by $w_i$ then $y=x_{i+1}$.
\end{itemize}
Then one sets $\Phi(\gamma)(x) = x_n$.

Let $r > 0$. A rooted colored $r$-ball around $x\in V$, $B_r(x)$ is a spanned subgraph
of $G$ with the induced $4$-coloring such that
\begin{itemize}
\item 
$d_G(x,y)\leq r$, for any $y\in B_r(x)$, where $d_G$ is the usual
shortest path metric.
\end{itemize}
Let $U^r_G$ be the set of rooted equivalence classes of $4$-colored
rooted $r$-balls in $G$. That is we consider two rooted $r$-balls
$B_r(x)$ and $B_r(y)$ equivalent if there exists a
graph isomorphism $\theta:B_r(x)\to B_r(y)$ preserving
the edge-colorings, mapping $x$ to $y$. Clearly $U^r_G$ is a finite set.
For each $x\in V$ the $r$-type of $x$ is the element $\alpha_r(x)\in U^r_G$ 
representing the $r$-ball around $x$.
Following Gromov \cite{Gro}, we call a $4$-colored graph $G$ a {\it graph
of dense holonomy} if, for any $\alpha\in U^r_G$, there exists
an integer $m_\alpha$ such that any ball of radius $m_\alpha$ in
$G$ contains at least one vertex $x$ with $r$-type $\alpha_r(x) = \alpha$. Note that
in the theory of Delone-sets such graphs are called {\it repetitive}
 \cite{LaPle}.
Also, we call a $4$-colored graph $G$ {\it generic} if for any $x\neq y\in V$
there exists $r>0$ such that the $r$-types of $x$ and $y$ are different.

\subsection{Amenable graphs}
Recall that an infinite, connected graph $G(V,E)$ of bounded vertex
degree is amenable if for any $\epsilon>0$, there exists a finite subset
$\Omega_\epsilon\subset V$ such that
$$\frac{|\partial \Omega_\epsilon|}{|\Omega_\epsilon|}<\epsilon\,,$$
where
$$\partial\Omega_\epsilon:=\{x\in\Omega_\epsilon\,\mid\,
\mbox{there exists $y\in V\backslash \Omega_\epsilon$
such that $x$ and $y$ are adjacent}\}$$ is the boundary of $\Omega_\varepsilon$.
We call a graph $G(V,E)$ {\it strongly amenable}
if there exists an increasing sequence of finite subsets
$\Omega_1 \subseteq
\Omega_2\subseteq\dots$ such that
\begin{itemize}
\item $\cup^\infty_{n=1} \Omega_n=V\,$ (exhaustion);
\item $\limn \frac{|\partial \Omega_n|}{|\Omega_n|}=0\,.$
\end{itemize}
Our first step shall be to construct strongly amenable generic
$4$-colored graphs with dense holonomy.

\subsection{The associated topological action}
For any graph $G(V,E)$ edge-colored by the set $\{A,B,C,D\}$ we have
a naturally associated compact metric space $X_G$ with a topological
$\bZ_2*\bZ_2*\bZ_2*
\bZ_2$-action by the construction below. Let
$\widetilde{\alpha}=\{\alpha_1\prec
\alpha_2\prec\dots\}$ be an infinite chain, such that $\alpha_r\in U^r_G$
and the $r$-ball around the root of $\alpha_{r+1}$ is just $\alpha_r$
(that is $\alpha_r\prec\alpha_{r+1})$. The chains above form the
compact metric space $X_G$, the {\it type space}, where
$$d_{X_G}(\{\alpha_1\prec
\alpha_2\prec\dots\}, \{\beta_1\prec
\beta_2\prec\dots\})=2^{-r}\,$$
where $r$ is the smallest integer for which $\alpha_r\neq \beta_r$. We shall
show (Proposition \ref{p13}) that the $\Gamma$-action on $G$ extends to $X_G$ in a natural way. For
each $x\in V$, $\pi(x)=\{\alpha_1(x)\prec
\alpha_2(x)\prec\dots\}$ is the associated element in $X_G$, where
$\alpha_r(x)$ is the $r$-type of $x$. Obviously, $G$ is generic if and only
if $\pi:G\to X_G$ is an injective map. We prove (Proposition \ref{p15}) that if $G$ is both
generic and of dense holonomy, then $X_G$ is homeomorphic to the
Cantor-set and the associated $\Gamma$-action is minimal. We shall also
prove (Section \ref{s8}) that for the strongly amenable generic graph we construct, this
$\Gamma$-action has at least two different
 ergodic invariant measures. Finally (Theorem \ref{tmt}), we shall
show that for two of these ergodic invariant measures $\mu_1$ and $\mu_2$,
$c(\mu_1)\neq c(\mu_2)$, where $c(\mu_i)$ denotes the {\it cost} of the
measurable equivalence relation $(\Gamma,X_G,\mu_i)$. This shows that
$(\Gamma,X_G,\mu_1)$ and $(\Gamma,X_G,\mu_2)$ are not
orbit equivalent relations. It is important to note that
although $(\Gamma,X_G,\mu_i)$ are not hyperfinite equivalence relations,
all their leaves are amenable graphs. The existence of such measurable
equivalence relations was first observed by Kaimanovich \cite{Kai}.
\section{Amenable actions on countable sets} \label{sect2}
Recall the notion of amenable actions on countable sets. Let $\Gamma$ be
a countable group and $X$ be a countable set. 
A homomorphism $\Phi:\Gamma\to
S(X)$ from $\Gamma$ into the full permutation group of $X$ is called an action.
The action is amenable \cite{GM},\cite{GN} if there exists an invariant,
finitely additive probability measure on $X$, that is a (positive real valued) function $\mu$ on the subsets
of $X$ such that:
\begin{itemize}
\item $\mu(X)=1,$
\item $\mu(A\cup B)=\mu(A)+\mu(B)$ if $A\cap B = \emptyset,$
\item $\mu(\Phi(\gamma)(A))=\mu(A)$, for any $\gamma\in\Gamma$, 
$A\subseteq X$.
\end{itemize}
It is well-known that $\Phi$ is amenable if and 
only if the following \Fo-condition
holds:

\noindent
For any finite set $1\in K\subset \Gamma$ and $\e>0$ there exists
a finite subset $\Omega_{K,\e}\subset X$ such that
$$\frac{|\Phi(K)\Omega_{K,\e}|}{|\Omega_{K,\e}|}< 1+\e\,.$$
Now let $\Gamma$ be a finitely generated group with symmetric generating
system $S$. The Schreier-graph ${\mathcal G}_\Phi = {\mathcal G}(\Gamma, \Phi, X, S)$ of an action $\Phi$ is given as follows:
\begin{itemize}
\item $V({\mathcal G}_\Phi)=X\,.$
\item $(x,y)\in E({\mathcal G}_\Phi)$ if $x=\Phi(s)(y)$ for some $s\in S$.
\end{itemize}
Clearly, $\Phi$ is amenable if and only if ${\mathcal G}_\Phi$ is an amenable 
graph.
If $\Gamma$ is amenable, than any action of $\Gamma$ is amenable as well. 
On the other hand, several non-amenable groups have faithful, transitive
amenable actions on countable sets (see \cite{GN} and \cite{GM}).
In this paper we consider amenable $\bZ_2*\bZ_2*\bZ_2*
\bZ_2$-actions on amenable graphs with proper $4$-colorings.
From now on $\Gamma$ shall always denote the group $\bZ_2*\bZ_2*\bZ_2*
\bZ_2$. As we have seen in the Introduction for any $\Gamma$-action $\Phi$
there exists a natural associated compact metric space of chains $X_G$.
\begin{proposition}
\label{p13}
The action $\Phi$ extends to a topological action $\tilde{\Phi}$ 
on the metric space $X_G$.\end{proposition}
\proof First note that, by definition, the image of $\pi:X \to
     X_G$, where, for $x \in X$, $\pi(x) = \{\alpha_1(x) \prec \alpha_2(x) \prec \ldots\}$ 
     is dense in $X_G$. Let $\widetilde{\alpha} = \{\alpha_1 \prec \alpha_2 \prec \ldots\}
     \in X_G$ and consider a word $\underline{w} = w_n w_{n-1}\cdots w_1 \in \Gamma$. For all
     $i = 1,2,\ldots$ let $x_i \in X$ be such that $\alpha_{i+n}(x_i) = \alpha_{i+n}$. If
     $y_i = \Phi(\underline{w})(x_i)$ setting $\beta_i = \alpha_i(y_i)$ one has $\beta_i
     \prec \beta_{i+1}$ and one defines
$$\widetilde{\Phi}(\uw)(\widetilde{\alpha}):=
\{\beta_1\prec\beta_2\prec\dots\}\,.$$ Obviously $\widetilde{\Phi}(\uw)$
is a continuous map for all $\uw\in\Gamma$ and $\widetilde{\Phi}(\uw)
(\pi(x))=\pi(\Phi(\uw)(x))$\,.\qed
\begin{proposition}
\label{p15}
If $\Phi$ is an amenable action of $\Gamma$, then $\widetilde{\Phi}$
has invariant measure on $X_G$.
\end{proposition}
\proof
First fix an ultrafilter $\omega$ on $\bN$ and consider the associated
ultralimit
$$\limo:l^\infty(\bN)\to\bR\,.$$
Let $\{\Omega_n\}^\infty_{n=1}$ be a \Fo-sequence in $X$ with respect to
$\Phi$. That is for any $\gamma\in\Gamma$,
$$\limn\frac{|\Phi(\gamma)\Omega_n\cup \Omega_n|}{|\Omega_n|}=1\,.$$
Let $\alpha\in U^r_G$ be an $r$-type of the Schreier-graph.
Denote by $\tau_n(\alpha)$ the number of vertices in the set $\Omega_n$
which are of $r$-type $\alpha$. Then let
$$\mu(\underline{\alpha}):=\limo\frac{\tau_n(\alpha)}{|\Omega_n|}\,,$$
where $\underline{\alpha}\subset X_G$ is the clopen set of chains 
$\{\delta_1\prec\delta_2\prec\dots\}$, with $\delta_r=\alpha$.
Now let \\ $\beta^1_{r+1}, \beta^2_{r+1},\dots\beta^k_{r+1}$ be the set of
$(r+1)$-types such that $\alpha \prec\beta^i_{r+1}$.
Clearly,
$$\tau_n(\alpha)=\sum^k_{i=1}\tau_n(\beta^i_{r+1})\,.$$
That is $
\mu(\underline{\alpha})=\sum^k_{i=1}\mu(\underline{\beta^i_{r+1}})\,,$
where $\underline{\beta^i_{r+1}}\subset X_G$ is the clopen set of chains
with $(r+1)$-type $\beta^i_{r+1}$. This shows immediately that $\mu$ defines
a Borel-measure on the compact metric space $X_G$.
In order to prove that $\mu$ is $\widetilde{\Phi}$-invariant it is enough to show that
$$\mu(\widetilde{\Phi}(w)(\underline{\alpha}))=\mu(\underline{\alpha})\,,$$
where $w$ is one of the generators of $\Gamma$, namely $A$, $B$, $C$ or $D$. 
Let $\{s_1,s_2,\dots,s_l\}\subset U_G^{r+1}$ be the set of $(r+1)$-types in
the Schreier graph of the action with the following property.
The $r$-ball around $\Phi(w)(x_i)$ is just $\alpha$, where
$x_i$ is the root of $s_i$. Observe that
$\widetilde{\Phi}(w)(\overline{\alpha})= \coprod^l_{i=1}\underline{s_i}$.
Hence we need to prove that
\begin{equation} \label{e17}
\mu(\underline{\alpha})=\sum^l_{i=1} \mu(\underline{s_i})\,.
\end{equation}
In order to check (\ref{e17}) it is enough to see that
\begin{equation}
\label{e18}
\limn\frac{|\tau_n(\alpha)-\sum_{i=1}^l \tau_n(s_i)|}{|\Omega_n|}=0\,.
\end{equation}
Let $\Omega_n(\alpha)$ be the set of vertices in $\Omega_n$ of $r$-type
$\alpha$ and $\Omega_n(s_i)$ be the set of vertices in 
$\Omega_n$ of $(r+1)$-type $s_i$. Let us observe that
\begin{itemize}
\item $\Omega_n(s_i)\cap \Omega_n(s_j)=\emptyset$ if $i\neq j$.
\item For any $x\in\Omega_n(\alpha)$ if $\Phi(w)(x)\in \Omega_n$, then
$\Phi(w)(x)\in\Omega_n(s_i)$ for some $i$.
\item For any $x\in \Omega_n(s_i)$  if
$\Phi(w)(x)\in \Omega_n$, then $\Phi(w)(x)\in\Omega_n(\alpha)$ (recall that
      $w$ is an involution, namely $w^2 = 1$). \end{itemize}

Therefore,
$$|\tau_n(\alpha)-\sum_{i=1}^l \tau_n(s_i)|\leq 2|\partial\Omega_n|\,.$$
Thus (\ref{e18}) immediately follows from the \Fo-property.\qed
\section{Generic actions of dense holonomy}
Again let $\Gamma$ be the free product of four cyclic groups of order two
$\bZ_2 * \bZ_2 * \bZ_2 * \bZ_2 $ and $\Phi$ be an amenable $\Gamma$-action
on the countable set $X$,
induced by a $4$-coloring of a graph $G(X,E)$. 
We call the action $\Phi$ generic,
resp. of dense holonomy, if the $4$-colored graph is generic, resp. of
dense holonomy. 
\begin{proposition}
\label{p20}
If $\Phi$ is a generic action on $X$ of dense holonomy, then $X_G$ is
homeomorphic to  a Cantor-set and the associated action $\overline{\Phi}$
is minimal. \end{proposition}
\proof
Let $\alpha\in U^r_G$, then there exists $\beta\neq \gamma\in U^s_G$
(for some $s>r$) such that $\alpha\prec\beta$ and $\alpha\prec\gamma$, 
that is the $r$-ball around the root of $\beta$ and $\gamma$ is of
type $\alpha$. Indeed, let $x\in X$ a vertex of $r$-type $\alpha$. Then by
the dense holonomy property there exists $y\in X$ of the same $r$-type. By
genericity there exists some $s>r$ such that the $s$-types of $x$ and $y$
are different. Consequently, $X_G$ is a compact metric space with
no isolated point, thus it is homeomorphic to the Cantor-set. 

\noindent
Now we prove minimality. Let $\widetilde{\alpha}=\{\alpha_1\prec
\alpha_2\prec\dots\}$ and $\widetilde{\beta}=\{\beta_1\prec
\beta_2\prec\dots\}$ be elements of $X_G$. 
 It is enough to prove that for any $r>0$ there exists
$\uw\in\Gamma$ such that $\widetilde{\Phi}(\uw)(\widetilde{\alpha})=
\{\beta_1\prec\beta_2\prec\dots\prec\beta_r\prec\gamma_{r+1}\dots\}$
that is $d(\widetilde{\beta}, \widetilde{\Phi}(\uw)(\widetilde{\alpha}))
\leq 2^{-r}\,.$ Let $m>r$ be an integer such that, by dense holonomy, any $m$-ball
in the graph $G$ contains a vertex of $r$-type $\beta_r$. 
Consider a vertex $p\in X$ such that its $2m$-type is $\alpha_{2m}$.
Let $q\in B_m(p)$ be a vertex of $r$-type $\beta_r$. Then there exists
a word $\uw=w_l w_{l-1}\dots w_1$, $l\leq m$
such that $\Phi(\uw)(p)=q\,.$ By the definition of the induced topological action $\tilde{\Phi}$ one has
$\widetilde{\Phi}(\uw)(\widetilde{\alpha})=
\{\beta_1\prec\beta_2\prec\dots\prec\beta_r\prec\dots\}$. \qed

\vskip 0.1in
\noindent
{\bf Note:} If the action of a group is generic but
not of dense holonomy the associated action on the type space might not be minimal.
Consider the group $\Delta=\bZ_2 * \bZ_2$, the dihedral group generated
by symbols $A$ and $B$. We define an action on the positive integers
the following way:

\begin{itemize}
\item $\Phi(A)(2n-1)=2n\,, n\geq 1\,.$
\item $\Phi(A)(2n)=2n-1\,, n \geq 1\,.$
\item $\Phi(B)(2n-1)=2n-2\,, n> 1\,.$
\item $\Phi(B)(2n)=2n+1\,, n> 1\,.$
\item $\Phi(B)(1)=1\,.$ \end{itemize}
It is easy to check that $\Phi$ is generic. The associated compact metric
space is 
$\{1,2,3,\dots\}\cup\{\infty\}$, where $$\widetilde{\Phi}(A)(\infty)=
\infty\quad \mbox{and} \quad \widetilde{\Phi}(B)(\infty)=
\infty\,.$$ Thus the action $\widetilde{\Phi}$ is not minimal.
\section{Technicalities}
For our main construction we need some simple lemmas.
\begin{lemma}
\label{l24}
Let $G(V,E)$ be a finite connected graph and $d$ be an integer
larger than $2$. Suppose that $A\subset V$ is a $2d$-net, that is if $x\neq
y\in A$ then $d_G(x,y)\geq 2d$. Also suppose that $diam\, G\geq 2d$.
Then $\frac{|A|}{|V|}\leq \frac{1}{d}\,.$
\end{lemma}
\proof
By the diameter condition, for each $x \in A$
   and $i=1,2,\ldots, d$ there exists $z_{x,i} \in V$ such that $d_G(x,z_{x,i})=i$. 
   Hence $|B_d(x)|\geq d$. Since the $d$-balls 
around the elements of a $2d$-net are disjoint, the lemma follows. \qed
\begin{lemma}
\label{l25}
Let $s_1<s_2<\dots<s_n$ be integers, $s_i\geq 10^{i+1}$. Suppose
that $diam(G)>10 s_n$. Then there exists a partition
$$V=R_1\coprod R_2\coprod\dots\coprod R_n\coprod R_{n+1}$$ such that
\begin{itemize}
\item For any $1\leq i \leq n$, $R_i$ is a $2s_i$-net.
\item For any $x\in V$, there exists $p\in R_i$ such that
$d_G(p,x)\leq 10 s_i$.
\end{itemize} \end{lemma}
\proof Let $R_1\subset V$ be a maximal $2s_1$-net in $V$. By the
previous lemma
$$\frac{|R_1|}{|V|}\leq \frac {1}{s_1}\,.$$
Moreover, for the same reason, for any spanned  subgraph $T\subseteq G$, 
with $diam(T)\geq 2s_1$:
$$\frac{|T\cap R_1|}{|T|}\leq \frac{1}{s_1}\,.$$
By maximality, for any $y\in V$ there exists $p\in R_1$ such that 
$d(p,y)\leq 10s_1$.
Now let $R_2\subset V$ be a $2s_2$-net which is maximal with respect to
the property that $R_1\cap R_2=\emptyset$. Again, for any spanned subgraph
$T\subseteq G$ with $diam(T)\geq 2s_2$:
\begin{equation}\label{e26}
\frac{|T\cap R_1|}{|T|}\leq \frac{1}{s_1}\quad,\quad
\frac{|T\cap R_2|}{|T|}\leq \frac{1}{s_2}\,.
\end{equation}
Let $x$ be an arbitrary vertex of $V$, we need to prove that
$B_{10s_2}(x)\cap R_2$ is non-empty. If $B_{10s_2}(x)\cap R_2=\emptyset$
then by maximality, $B_{5s_2}(x)$ is completely filled with vertices
from the set $R_1$. However, this is in contradiction with (\ref{e26}). 
Now we can construct the partition by a simple induction. \qed

\section{The main construction}

The goal of this section is to construct a
\begin{itemize}
\item faithful, transitive and amenable
\item generic
\item dense holonomy
\end{itemize}
$\Gamma$-action on a countable set $X$. In the whole section $\Gamma$
denotes the group $\bZ_2*\bZ_2*\bZ_2*
\bZ_2$ generated by the elements $\{A,B,C,D\}$ of order two, and
$\widetilde{\Gamma}$ denotes the subgroup generated by $\{A,B,C\}$.
We also fix a parameter $m>10$ to be specified later. We also need
two auxilliary graph sequences. The
first one is $\{C_i\}^\infty_{i=1}$, where $C_i$ is the cycle of
length $2i$ edge-colored by $A$ and $B$. The second graph
sequence is given the following the way.
The group $\widetilde{\Gamma}$ is residually finite (free product
of finite groups), hence one can pick a decrasing sequence of finite index normal subgroups
      with trivial intersection:
$$\widetilde{\Gamma}\rhd N_1 \rhd N_2 \rhd N_3\rhd\dots\quad
,\cap^\infty_{i=1} N_i=\{1\}\,.$$
Then the $3$-colored graph $K_i$ is the Cayley-graph of the finite
group $\widetilde{\Gamma}/N_i$. Note that $\{K_i\}^\infty_{i=1}$
is a large girth sequence, that is, for any $s > 0$, if $i$ is large enough, 
then $K_i$ does not contain cycles of length less than $s$. 
Finally, we lexicographically enumerate the words in $\Gamma$: $A, B, C, D, AB, AC,
AD, BA, \ldots$ and denote the $n$-th word by 
${\underline{w}}_n = w_{k_n}^{(n)}  w_{k_n-1}^{(n)}\cdots w_2^{(n)}w_1^{(n)}$; note that the
word-length $k_n$ of ${\underline{w}}_n$ is at most $n$.
     
\noindent
\underline{\bf {Step (1)}} \quad  Let $G_0(V_0,E_0)$ be the graph consisting of one single
vertex $p$. Let $G_1(V_1,E_1)$ be the graph constructed the following
way. $H_1$ is a triangle, edge-colored by $A,B$ and $C$. The vertex $p$
is connected to a vertex $q$ of $H_1$ by an edge colored by $D$.
Finally a vertex of $H_1$ which is not $q$ shall be called $r_1$.
Thus we have graphs $G_0\subset G_1$.
In the $n$-th step we shall have graphs $G_0\subset G_1\subset G_2\subset
\dots\subset G_n$. Finally let us set $s_1:=2m$, $s_2=m^2 2^{|V(G_1)|}$.

\noindent
\underline {\bf {Step (2)}} \quad Let $H_2$ be a graph $C_i$, in such a way that
$diam\,C_i>10 s_2$. By Lemma \ref{l25} there exists a partition
$$V(H_2)=R^2_1\coprod R^2_2\coprod R^2_3\,$$
where $R^2_1$ is a $2s_1$-net, $R^2_2$ is a $2s_2$-net satisfying the 
properties described in the lemma. Connect each point of $R^2_1$ to
a copy of $G_0$ (a single vertex) with an edge colored by $D$. For one single vertex
of $R^2_2$ connect a copy of $G_1$ (at the vertex $r_1$ !) with an edge
colored by $D$. Let $x_2$ be a vertex of $R^2_3$. Consider the first word $\underline{w}_1 = A$. 
Let $y_0^2, y_1^2$ be two new vertices (i.e.
distinct from any other vertex involved in the construction till now) and
connect $x_2$ to $y_0^2$ by an edge colored by $D$ and connect $y_0^2$ to $y_1^2$ by an
edge colored by $\underline{w}_1 = A$. Finally pick an other vertex of
$R^2_3$ and call it $r_2$.

\noindent
Thus we constructed $4$-colored graphs $G_0\subset G_1\subset G_2$.
In $G_2$ we have a spanned subgraph $H_2$ such that
$$\frac{|V(G_2)\backslash V(H_2)|}{|V(G_2)|}<\frac{1}{m}\,.$$
Set $s_3:=m^3 2^{|V(G_2)|}$\,.

Indeed, $|V(G_2)| = 2 |R_1^2| + |R_2^2| + |R_3^2| + 2$ (the last term corresponds to 
      $|\{y_0^2, y_1^2\}|$) and $$\frac{|V(G_2)\setminus V(H_2)|}{|V(G_2)|} =
      \frac{|R_1^2| + 2}{|V(G_2)|} \leq \frac{|R_1^2| + 2}{|V(H_2)|} \leq 
      \frac{1}{2m} + \frac{1}{10m^2 2^4} < 1/m.$$

\noindent
\underline{\bf {Step (n)}}\quad Suppose that we have already constructed
the graphs $G_0\subset G_1\subset G_2\subset\dots\subset G_{n-1}\,$ with
the following properties:
\begin{itemize}
\item
If $i>1$, then $G_i$ contains a spanned subgraph $H_i$ such that
$\frac{|V(G_i)\backslash V(H_i)|}{|V(G_i)|}<\frac{1}{m}\,.$
\item  $s_i = m^i 2^{|V(G_{i-1})|}$.
\item For each $H_i$, we fixed a vertex $r_i$.
\end{itemize}
Now if $n$ is even, then let $H_n$ be isomorphic to $C_i$ for some $i$,
if $n$ is odd, then let $H_n$ be isomorphic to $K_i$ for some $i$,
satisfying the inequality $diam (H_n)\geq 10 s_{n}$.
We apply Lemma \ref{l25} again.
$$V(H_n)=R^n_1\coprod R^n_2\coprod\dots\coprod R^n_{n+1}\,,$$
where $R^n_i$ is a $2s_i$-net and for each $x\in V(H_n)$ there exists
$p\in R^n_i$ such that
$$d_{H_n}(p,x)<10s_i\,.$$
For $i=1,2,\ldots,n-1$ we connect a copy of $G_{i-1}$ (at vertex $r_i$)
to each vertex of $R^n_i$ by an edge
colored by $D$. Then we connect one single copy of $G_{n-1}$ to $H_n$ by an edge colored by $D$ between
$r_{n-1}$ and an arbitrary vertex in $R_n^n$.
Now we pick a vertex $x_n$ of $R^n_{n+1}$. Consider the $n$th word
$\underline{w} = w^{(n)}_{k_n}w^{(n)}_{k_n-1}\ldots w^{(n)}_1$.
\begin{itemize}\item
If $w^{(n)}_{1}=D$, then choose vertices $y_1^{n},
y_2^{n},\dots,y^{n}_{k_{n}}$ disjoint
 from any previous vertices and connect the
path $(y^{n}_1,y^{n}_2,\dots,y^{n}_{k_n})$ to
  $x_n$ by the edge $(x_n,y^n_1)$ 
colored by $D$ and let color the edge $(y^n_i,y^n_{i+1})$ by $w^{n}_{i+1}$.
\item
If $w^{(n)}_{1}\neq D$, then choose vertices $y_0^{n},y_1^{n},
y_2^{n},\dots,y^{n}_{k_{n}}$ disjoint
 from any previous vertices and connect the
path $(y^{n}_0,y^{n}_1,\dots,y^{n}_{k_n})$ to
  $x_n$ by the edge $(x_n,y^n_1)$ 
colored by $D$ and let color the edge $(y^n_i,y^n_{i+1})$ by $w^{n}_{i+1}$.
\end{itemize}
Finally pick an other vertex $r_n$ form $R^n_{n+1}\,.$

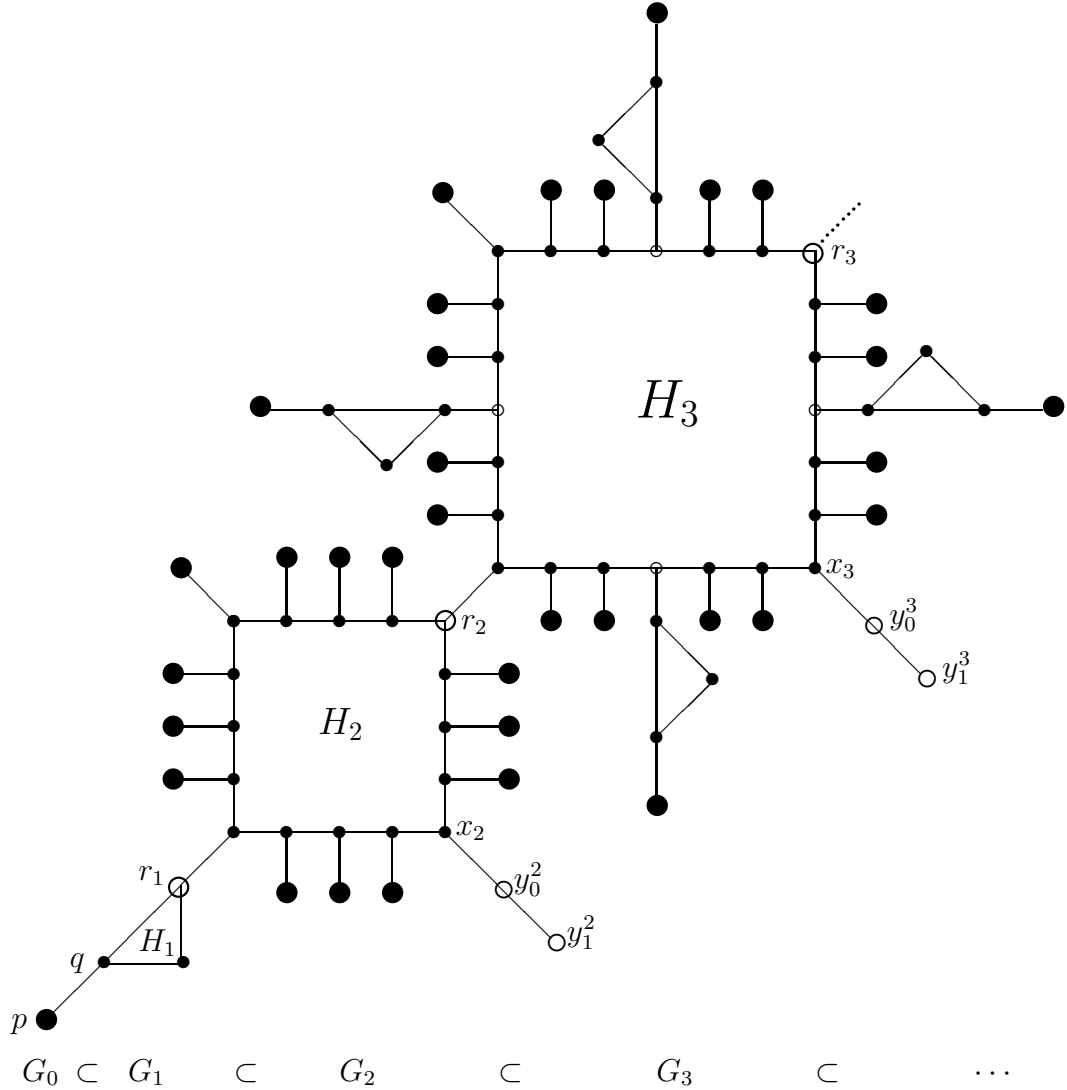
\begin{center}
\begin{picture}(480,400)
\put(30,30){\line(1,1){70}}
\put(50,50){\line(1,0){30}}
\put(80,50){\line(0,1){30}}
\put(100,100){\line(1,0){80}}
\put(100,100){\line(0,1){80}}
\put(100,180){\line(1,0){80}}
\put(180,100){\line(0,1){80}}
\put(180,180){\line(1,1){20}}
\put(200,200){\line(1,0){120}}
\put(200,200){\line(0,1){120}}
\put(200,320){\line(1,0){120}}
\put(320,200){\line(0,1){120}}

\put(120,100){\line(0,-1){20}}
\put(140,100){\line(0,-1){20}}
\put(160,100){\line(0,-1){20}}
 
\put(120,180){\line(0,1){20}}
\put(140,180){\line(0,1){20}}
\put(160,180){\line(0,1){20}}

\put(100,120){\line(-1,0){20}}
\put(100,140){\line(-1,0){20}}
\put(100,160){\line(-1,0){20}}

\put(180,120){\line(1,0){20}}
\put(180,140){\line(1,0){20}}
\put(180,160){\line(1,0){20}}

\put(100,180){\line(-1,1){20}}

\put(180,100){\line(1,-1){40}}

\put(220,200){\LARGE \line(0,-1){20}}
\put(240,200){\line(0,-1){20}}
\put(260,200){\line(0,-1){86}}
\put(260,136){\line(1,1){22}}
\put(260,180){\line(1,-1){22}}
\put(280,200){\line(0,-1){20}}
\put(300,200){\line(0,-1){20}}

\put(220,320){\line(0,1){20}}
\put(240,320){\line(0,1){20}}
\put(260,320){\line(0,1){86}}
\put(260,340){\line(-1,1){22}}
\put(260,384){\line(-1,-1){22}}
\put(280,320){\line(0,1){20}}
\put(300,320){\line(0,1){20}}

\put(200,220){\line(-1,0){20}}
\put(200,240){\line(-1,0){20}}
\put(200,260){\line(-1,0){86}}
\put(180,260){\line(-1,-1){22}}
\put(136,260){\line(1,-1){22}}
\put(200,280){\line(-1,0){20}}
\put(200,300){\line(-1,0){20}}

\put(320,220){\line(1,0){20}}
\put(320,240){\line(1,0){20}}
\put(320,260){\line(1,0){86}}
\put(340,260){\line(1,1){22}}
\put(384,260){\line(-1,1){22}}
\put(320,280){\line(1,0){20}}
\put(320,300){\line(1,0){20}}

\put(200,320){\line(-1,1){20}}
\put(320,200){\line(1,-1){40}}

\put(24,24){\LARGE $\bullet$}

\put(155,72){\LARGE $\bullet$}
\put(115,72){\LARGE $\bullet$}
\put(135,72){\LARGE $\bullet$}

\put(155,199){\LARGE $\bullet$}
\put(115,199){\LARGE $\bullet$}
\put(135,199){\LARGE $\bullet$}

\put(72,115){\LARGE $\bullet$}
\put(72,135){\LARGE $\bullet$}
\put(72,155){\LARGE $\bullet$}

\put(199,115){\LARGE $\bullet$}
\put(199,135){\LARGE $\bullet$}
\put(199,155){\LARGE $\bullet$}

\put(215,175){\LARGE $\bullet$}
\put(235,175){\LARGE $\bullet$}
\put(275,175){\LARGE $\bullet$}
\put(295,175){\LARGE $\bullet$}

\put(215,338){\LARGE $\bullet$}
\put(235,338){\LARGE $\bullet$}
\put(275,338){\LARGE $\bullet$}
\put(295,338){\LARGE $\bullet$}

\put(172,215){\LARGE $\bullet$}
\put(172,235){\LARGE $\bullet$}
\put(172,275){\LARGE $\bullet$}
\put(172,295){\LARGE $\bullet$}

\put(338,215){\LARGE $\bullet$}
\put(338,235){\LARGE $\bullet$}
\put(338,275){\LARGE $\bullet$}
\put(338,295){\LARGE $\bullet$}

\put(75,195){\LARGE $\bullet$}

\put(48,48){$\bullet$}
\put(78,48){$\bullet$}
\put(74,74){\LARGE $\circ$}

\put(97,97){$\bullet$}
\put(117,97){$\bullet$}
\put(137,97){$\bullet$}
\put(157,97){$\bullet$}
\put(177,97){$\bullet$}

\put(97,177){$\bullet$}
\put(117,177){$\bullet$}
\put(137,177){$\bullet$}
\put(157,177){$\bullet$}
 
\put(97,117){$\bullet$}
\put(97,137){$\bullet$}
\put(97,157){$\bullet$}
\put(97,177){$\bullet$}

\put(177,117){$\bullet$}
\put(177,137){$\bullet$}
\put(177,157){$\bullet$}
\put(175,175){\LARGE $\circ$}

\put(197,197){$\bullet$}
\put(217,197){$\bullet$}
\put(237,197){$\bullet$}
\put(257,197){$\circ$}
\put(277,197){$\bullet$}
\put(297,197){$\bullet$}
\put(317,197){$\bullet$}

\put(197,317){$\bullet$}
\put(217,317){$\bullet$}
\put(237,317){$\bullet$}
\put(257,317){$\circ$}
\put(277,317){$\bullet$}
\put(297,317){$\bullet$}
\put(314,314){\LARGE $\circ$}

\put(321,321){$\cdot$}
\put(323,323){$\cdot$}
 
\put(325,325){$\cdot$}

\put(327,327){$\cdot$}
\put(329,329){$\cdot$}
\put(331,331){$\cdot$}

\put(333,333){$\cdot$}
\put(335,335){$\cdot$}

\put(218,54){\Large$\circ$}
\put(198,74){\Large$\circ$}

\put(358,154){\Large$\circ$}
\put(338,174){\Large$\circ$}

\put(64,55){$H_1$}

\put(132,137){\large $H_2$}

\put(252,257){\LARGE $H_3$}
\put(16,26){$p$}
\put(38,49){$q$}

\put(64,82){$r_1$}

\put(184,99){$x_2$}
\put(186,177){$r_2$}
\put(206,79){$y_0^2$}
\put(226,59){$y_1^2$}

\put(326,317){$r_3$}

\put(324,198){$x_3$}
\put(348,179){$y_0^3$}
\put(368,159){$y_1^3$}

\put(197,217){$\bullet$}
\put(197,237){$\bullet$}
\put(197,257){$\circ$}
\put(197,277){$\bullet$}
\put(197,297){$\bullet$}

\put(317,217){$\bullet$}
\put(317,237){$\bullet$}
\put(317,257){$\circ$}
\put(317,277){$\bullet$}
\put(317,297){$\bullet$}

\put(257,177){$\bullet$}
\put(257,133){$\bullet$}
\put(278,155){$\bullet$}
\put(255,105){\LARGE $\bullet$}

\put(257,337){$\bullet$}
\put(257,381){$\bullet$}
\put(235,359){$\bullet$}
\put(255,405){\LARGE $\bullet$}

\put(177,257){$\bullet$}
\put(133,257){$\bullet$}
\put(155,236){$\bullet$}
\put(105,256){\LARGE $\bullet$}

\put(337,257){$\bullet$}
\put(381,257){$\bullet$}
\put(359,279){$\bullet$}
\put(405,256){\LARGE $\bullet$}

\put(174,337){\LARGE $\bullet$}

\put(20,6){$G_0$}
\put(60,6){$G_1$}
\put(140,6){$G_2$}
\put(260,6){$G_3$}

\put(40,6){$\subset$}
\put(100,6){$\subset$}
\put(200,6){$\subset$}
\put(320,6){$\subset$}
\put(380,6){$\cdots$}

\end{picture}
\end{center}

Thus we constructed a sequence of graphs
$G_0\subset G_1\subset \dots\subset G_n$. We proceed by induction.
Denote the $4$-colored graph $\cup^\infty_{n=1} G_n$ by $G$. 
In the next section we prove that $G$ is generic, has dense holonomy and carries a faithful, transitive, amenable action of the group
$\Gamma$.

\section{The graph $G$}
In this section we check that the graph $G$ constructed in the previous
section has indeed all the required properties.
\begin{proposition}
\label{p35}
The action of $\Gamma$ is faithful, transitive and (strongly) amenable.
\end{proposition}
\proof
Transitivity is obvious from the connectivity of $G$. Faithfulness was
ensured by the attachment of paths $(y^n_1,y^n_2,\dots,y^n_{k_n})$.
The subgraphs $\{G_n\}^\infty_{n=1}$ are forming an increasing and exhausting \Fo-sequence, since
$\partial G_n=\{r_n\}$. \qed
\begin{proposition}
\label{p35b}
The $4$-colored graph $G$ is generic. \end{proposition}
\proof First we need a lemma.
\begin{lemma}\label{l36}
The single vertex $p$ of $G_0$ has the property that for any $q\in V(G)$, $q \neq p$,
there exists $r$ such that $B_r(p)$ and $B_r(q)$ are not isomorphic (as colored graphs,
in other words $\alpha_r(p) \neq \alpha_r(q)$).
\end{lemma}
\proof
Let $q\in V(G)$ be a vertex of degree one (otherwise the statement trivially
holds). If $q$ was defined in {\bf Step (k)}, then there exists a spanned
connected subgraph of $G$ which is isomorphic to $H_k$ such that the unique
 path
from $q$ to the subgraph does not contain any vertex of a subgraph isomorphic
to $H_{k-1}$. Indeed, $H_k$ itself can be chosen as the spanned subgraph.
On the other hand any path connecting $p$ to a spanned subgraph isomorphic
to $H_k$ passes through such a vertex. This proves our lemma.
\qed

\noindent
Now we finish the proof of Proposition \ref{p35b}. Let $x, y$ be two distinct vertices of $G$. By transitivity of the action there exists $\overline{w} \in \Gamma$ such that $\Phi(\overline{w})(x) = p$ (the single vertex of $G_0$). Set $q = \Phi(\overline{w})(y)$. By the previous lemma, there exists $r$ such that $B_r(p)$ is not colored-isomorphic to $B_r(q)$.
Then $B_{r+|\uw|}(x)$ is not colored-isomorphic to $B_{r+|\uw|}(y)$, where
$|\uw|$ is the word-length of $\uw$. \qed
\begin{proposition}\label{p37}
The $4$-colored graph $G$ is of dense holonomy.
\end{proposition}
\proof
Let $t\in V(G)$ and $r>0$. Suppose that $B_r(t)\subset G_i$. Let $M=
diam (G_{i+3})\,.$ 
\begin{lemma} For any $x \in V(G)$, the ball $B_M(x)$ contains a vertex
     $y$ belonging to a spanned subgraph $L$ isomorphic to $H_j$ for some $j \geq i+3$.
\end{lemma}
\proof
If $x\in V(G_{i+2})$ then the statement obviously holds.
We proceed by induction.
Suppose that the statement holds for any $x\in V(G_l)$. Let 
$x\in V(G_{l+1})\backslash V(G_l)$. If $x\in H_{l+1}$ then
$L$ can be chosen as $H_{l+1}$ itself. If $x$ is attached to $H_{l+1}$
then we have two cases:
\begin{itemize}
\item $B_M(x)$ intersects $H_{l+1}$: then the statement is clearly true.
\item $B_M(x)$ does not intersect $H_{l+1}$,: then there exists
a vertex $x'\in G_l$ such that $B_M(x)$ and $B_M(x')$ are isomorphic, hence
the statement follows by induction. \end{itemize} 
Now if $y\in L$, $L\simeq H_j$, $j\geq i+3$. Then in a $(10 s_{i+3})$-neighborhood of $y$ we attached a copy of $G_{i+3}$. Hence
$B_{M+10s_{i+3}}(y)$ contains a vertex of the same $r$-type as $t$.
Therefore, $B_{2M+10s_{i+3}}(x)$ contains a vertex of the same $r$-type
as the vertex $t$, in other words (recalling the
      definition of dense holonomy) if $\alpha_r(t) = \alpha$, then $m_\alpha =  2M+10s_{i+3}$, and the proposition follows.
\qed
\section{Graphed equivalence relations}
Let us recall the notion of Borel-equivalence relations and their
$L$-graphings from the book of Kechris and Miller (\cite{KM}, Section 17).
Let $Z$ be a standard Borel-space. Let $\{A_i\}^\infty_{i=1}$ and
$\{B_i\}^\infty_{i=1}$ be two families of Borel subsets of $Z$ and $\phi_i:A_i\to B_i$ be partial
Borel-isomorphisms. Then the system $\{\phi_i\}^\infty_{i=1}$ defines
a countable Borel-equivalence relation $E$ on $Z$ in the following way: for $x,y \in Z$ we set $x \equiv_E y$ if there exists a finite
      sequence $x_1, x_2, \ldots, x_n \in Z$ such that 
\begin{itemize}
\item $x_1=x, x_n=y$ and
\item for $i=1,2,\ldots,n$ one has $x_{i+1} = \phi_j^{\pm 1}(x_i)$ for some $j \in {\mathbb N}$.
\end{itemize}
This also means that for each $z\in Z$ the system $\{\phi_i\}^\infty_{i=1}$
defines a $L$-graphing $G$, that is a graph structure on the equivalence
classes on $E$. The simplest examples of Borel-equivalence relations
are topological actions of finitely generated groups on compact metric spaces: 
for $x, y \in Z$,  $x \equiv_E y$ if there exists $\gamma \in \Gamma$ such that
     $x = \gamma(y)$. Moreover, if the group $\Gamma$ is finitely generated by a symmetric generating 
     system $S$, then $x$ is connected to $z$ if $x = s(z)$ for some $s \in S$.

\noindent
Let $\mu$ be a  Borel probability measure on $Z$ preserved by
 the maps $\{\phi_i\}^\infty_{i=1}$ e.g. an invariant measure of a
 $\Gamma$-action. Then we call $\mu$ an $E$-invariant measure. Note that
$E$-invariance does not depend on the choice of the $L$-graphing $\cG$
defining $E$. A {\it measurable equivalence relation} is given by
a triple $(Z,E,\mu)$, where $Z$ is a standard space, $E$ is a countable
Borel-equivalence relation and $\mu$ is an $E$-invariant measure \cite{Gab}.
The measurable equivalence relations $(Z_1,E_1,\mu_1)$ and $(Z_2,E_2,\mu_2)$
are orbit-equivalent if there exists a measure preserving map $T:Z_1\to Z_2$
such that for almost all $z\in Z_1$ the image of the class of $z$ is just the 
equivalence class of $T(z)$. 
Let $\phi:Z\to Z$ be a single homeomorphism with an ergodic measure $\mu$.
The orbit equivalence class of the relation $E_\phi$ is called the {\it hyperfinite}
ergodic relation. The ergodic actions of amenable groups always define
the hyperfinite relation. In \cite{GN}, the authors constructed minimal
actions of certain non-amenable groups with ergodic measure still defining
the hyperfinite relation. 
For an $L$-graphing $\cG$ and an $E$-invariant measure $\mu$ the {\it edge measure} of ${\mathcal G}$ is defined as
$$e({\mathcal G},\mu) = \frac{1}{2}\int_Z \deg_{\mathcal G}(x)d\mu(x),$$
where $\deg_{\cG} (x)$ is the degree of $x$ in the $L$-graphing.
The {\it cost} of the equivalence relation is defined as
$$c_\mu(E):=\inf_{\cG} e(\cG,\mu)\,,$$
where the infimum is taken over all $L$-graphings defining $E$. The cost
of the hyperfinite ergodic relation  is $1$. Obviously, orbit equivalent
relations have the same cost.
Suppose that $\wp$ is the topological action of the
finitely generated group $\Gamma$ on the compact metric space $X_G$,
preserving the measure $\mu$. Then the $L$-graphing of $\wp$
is just the union of the Schreier-graphs of the $\Gamma$-orbits.
Kaimanovich \cite{Kai} proved that if $\mu$ is an ergodic invariant
measure for the action then the hyperfiniteness of the induced relation
implies that all the leaves (the Schreier-graphs) are amenable.
\begin{proposition}
\label{p45}
Let $(\widetilde{\Phi},\Gamma,X_G)$ be the associated minimal action
given in Section \ref{sect2} using or main construction. Then all the
leaves are amenable.
\end{proposition}
\proof
First note that the action $\widetilde{\Phi}$ is regular on $X_G$, that is
if $p\in X_G$, $\gamma\in\Gamma$ then either $\widetilde{\Phi}
(\gamma)(p)\neq p$
or the homeomorphism $\widetilde{\Phi}
(\gamma)$ fixes a whole neighborhood of $p$. Indeed let
$p=\{\alpha_1,\prec \alpha_2\prec\dots\prec\alpha_s\prec\dots\}$ and
let the word-length of $\gamma$ be $s$.
If $\widetilde{\Phi}
(\gamma)(p)= p$, then 
$\widetilde{\Phi}
(\gamma)$ also fixes all $q \in X_G$ of the form $q =\{\alpha_1,\prec
\alpha_2\prec\dots\prec\alpha_s\prec\beta_{s+1}\prec\dots\}$.

\noindent
We say that two $4$-colored graphs $G_1,G_2$ are locally isomorphic
if $U^r_{G_1}=U^r_{G_2}$ for any $r>0$. Clearly if $G_1$ is amenable and of
dense holonomy and $G_2$ is locally isomorphic to $G_1$, then $G_2$ is
amenable and of dense holonomy as well \cite{Gro}. Since our graph $G$ is
generic, the dense leaf $\pi(G) \subset X_G$ is isomorphic to $G$. By minimality, regularity 
and dense holonomy, for any two orbits the corresponding Schreier-graphs are locally
isomorphic (see also \cite{GN}) and thus our proposition follows.\qed

\noindent
As we shall see later, our $\Gamma$-action $\tilde{\Phi}$ defines a non-hyperfinite
measurable equivalence relation with respect to some ergodic invariant
measure. The existence of non-hyperfinite relations with amenable leaves
was first observed by Kaimanovich \cite{Kai}.
\section{The invariant measures}\label{s8} 
Notice that in our main construction we have a parameter $m$ and we have in fact
constructed a graph $G\mi$ for each $m>10$.
By Lemma \ref{l24},
$$V(H\mi_n)=R^{n,(m)}_1\coprod R^{n,(m)}_2\coprod\dots\coprod R^{n,(m)}_{n+1}\,,$$
where
$$\frac{|R^{n,(m)}_i|}{|V(H\mi_n)|} < \frac{1}{m^i 2^{|V(G_{i-1}^{(m)})|}}$$
for $i = 1,2,\ldots, n$ so that
$$\frac{|R^{n,(m)}_{n+1}|}{|V(H\mi_n)|}>1-\frac{1}{m}\,.$$
Here the upper index $(m)$ indicates that the parameter is chosen to be $m$.
For each vertex of $R^n_i$ we attached a copy of a graph of size $|V(G_{i-1}^{(m)})|$\,.
We also attached a path of length at most $n$ to one vertex of
$R^{n,(m)}_{n+1}$. Thus we have the following inequality:
\begin{equation}\label{e47}
\frac{|V(G\mi_n)\backslash V(H\mi_n)|}{|V(H\mi_n)|}<\frac{1}{m}\,.
\end{equation}
Indeed
    $|V(G_n^{(m)})| \leq 2|R_1^{n,(m)}| + |V(G_1^{(m)})|\cdot |R_2^{n,(m)}|+ \cdots +
    |V(G_{n-1}^{(m)})|\cdot |R_n^{n,(m)}|+ n$ where the last term majorizes $k_n$, the length
     of the $n-$th word $\underline{w}$. Thus,
    $$\frac{|V(G_n^{(m)}) \setminus V(H_n^{(m)})|}{|V(H_n^{(m)})|} \leq
    \sum_{i=1}^{n}  \frac{|V(G_{i-1}^{(m)})| -1}{m^i 2^{|V(G_{i-1}^{(m)})|}} < \frac{1}{m}.$$
    
Let $\Omega\mi_n:=V(G\mi_n)\backslash V(G\mi_{n-1})\,.$ Clearly
$|\partial 
\Omega\mi_n|=2$, thus $\{\Omega\mi_n\}^\infty_{n=1}$ forms a \Fo-sequence
in $G\mi$.
We consider two \Fo-sequences for $G\mi$: $\{F^{1,(m)}_n\}^\infty_{n=1}$
and $\{F^{2,(m)}_n\}^\infty_{n=1}$, where $F^{1,(m)}_n=\Omega\mi_{2n}$
and $F^{2,(m)}_n=\Omega\mi_{2n+1}$.
By the averaging process described in Section \ref{sect2} we then obtain
two invariant measures on $X_G$, $\mu_1^{(m)}$ and $\mu_2^{(m)}$, respectively. 
Thus, for each $m>10$ we have
\begin{itemize}
\item A compact metric space $X_{G\mi}$.
\item A minimal $\Gamma$-action $\widetilde{\Phi}\mi$.
\item An $L$-graphing $\cG\mi$ associated with the action.
\item A Borel-equivalence relation $E\mi$ defined by $\cG\mi$.
\item Two $E\mi$-invariant measures: $\mu\mi_1$ and $\mu\mi_2$.
\end{itemize}
\begin{proposition}\label{p51}
$\lim_{m\to\infty} e(\cG\mi,\mu\mi_1)=1\,.$
\end{proposition}
\proof
Let $Y\mi\subset X_{\cG\mi}$ be the set of points $x$ such that
$\deg_{\cG\mi}(x)>2$. It is enough to prove that
$\lim_{m\to\infty} \mu\mi_1(Y\mi)=0\,.$
By the definition of the probability measure $\mu\mi_1$,
$$\mu\mi_1(Y\mi)=\limo \frac{|L\mi_n|}{|F^{1,(m)}|}\,,$$
where $L\mi_n$ is the set of vertices in $\Omega\mi_{2n}$ having degree
larger than $2$. By (\ref{e47}):
$$|L\mi_n|\leq\frac{2}{m} |V(H\mi_{2n})|\leq \frac{2}{m}
 |F^{1,(m)}|\,.\quad\qed$$
\section{Lower estimate for the cost}
A point $x\in X_{G\mi}$ is called $\widetilde{\Gamma}$-free if for any
$\gamma\in \widetilde{\Gamma} \setminus \{1\}$, $\widetilde{\Phi}(\gamma)(x)\neq x\,.$
Clearly the set of $\widetilde{\Gamma}$-free points $A\mi\subset X_{G\mi}$
is closed.
\begin{lemma} \label{l54}
$\lim_{m\to\infty} \mu\mi_2(A\mi)=1\,.$
\end{lemma}
\proof
Let $A\mi_k$ be the set of points in  $X_{G\mi}$ for which
$\widetilde{\Phi}(\gamma)(x)\neq x$ whenever the word-length of $\gamma
\in \widetilde{\Gamma}$
is not greater than $k$. Clearly, $\cap^\infty_{k=1} A\mi_k=A\mi\,.$
Observe that
$$\mu\mi_2(A\mi_k)=\limo\frac{|W\mi_{n,k}|}{|F^{2,(m)}_n|}\,,$$
where $W\mi_{n,k}$ is the set of vertices $p$ in $F^{2,(m)}_n$
such that $\Phi(\gamma)(p)\neq p$, if $|\gamma|\leq k$.
Obviously, if $p\in H\mi_n$ and $n$ is a large odd number, then
$p\in W\mi_{n,k}$. Thus $\mu\mi_2(A\mi_k)\geq 1-\frac{2}{m}\,.$
Therefore $\mu\mi_2(A\mi)\geq 1-\frac{2}{m}$ as well. \qed

\noindent
Let $E$ be a Borel-equivalence relation on a compact metric space $Z$
and $\mu$ be an $E$-invariant measure. Let $Y\subseteq Z$ be a Borel-set
of positive measure. Then $(E\mid Y)$ is the induced equivalence relation
and $\mu_Y$ is the restriction of $\mu$ onto $Y$. We call $Y$ a
{\it complete} Borel-section if $Y$ intersects almost every classes of $E$.
By a lemma of Gaboriau (\cite{KM},  Theorem 21.1) if $Y$ is a complete section then
\begin{equation}\label{e55}
c_\mu(E)=c_{\mu_Y}(E\mid Y)+\mu(Z\backslash Y)
\end{equation}
\begin{lemma}
\label{l55}
The induced equivalence relation $(E\mid A\mi)$ on the set of 
$\widetilde{\Gamma}$-free points is exactly the relation given 
by the $\widetilde{\Gamma}$-action.
\end{lemma}
\proof
Notice that in the graph $G\mi$ there is no simple cycle containing
an edge colored by $D$. Hence the same holds for the orbits
in $X_{G\mi}$. Thus if $p,q\in A\mi$ and $\gamma$ contains the symbol
$D$, then $\widetilde{\Phi}(\gamma)(p)\neq q$.\qed
\begin{lemma} \label{l56}
$$c_{\mu_2^{(m)}}(E | A^{(m)}) =\frac{3}{2}\mu\mi_2(A\mi)\,.$$
\end{lemma}
\proof
The lemma is a simple consequence of the famous theorem of Gaboriau \cite{Gab}
(see also  \cite{KM}) stating that the cost of an equivalence
relation defined by an  $L$-treeing $\cG$ is exactly the edge measure
of $\cG$. Simply note that the edge measure of a free $Z_2*Z_2*Z_2$- action on
a compact metric space with an invariant measure is $\frac{3}{2}$ times
the measure of the space. \qed

\begin{lemma} \label{l57}
$A\mi$ is a complete section of $E$.
\end{lemma}
\proof
We need to prove that the set of points $p\in X_{G\mi}$, such that
there exists no $\gamma\in\Gamma$ with $\widetilde{\Phi}(\gamma)(p)\in A\mi$
has $\mu\mi_2$-measure zero. We denote this set by $B\mi$.

\noindent
Let $S_{k,n}\subset F^{2,(m)}_n$ be the set of vertices $p$ such that
$d_{G_n}(p,H\mi_{2n+1})\leq k$.
By our construction of $G\mi$ it is clear that for any $\epsilon>0$
there exists $k_\epsilon$ such that
$$\frac{|S_{k_\epsilon,n}|}{|F^{2,(m)}_n|}\geq 1-\epsilon\,.$$
Therefore $\mu\mi_2(\Theta_{k_\epsilon})\geq 1-\epsilon$, where $\Theta_{k_\varepsilon}$ is the set of points $x \in X_{G^{(m)}}$
     for which there exists $\delta \in \Gamma$ with $|\delta| \leq k_\varepsilon$
     such that $\tilde{\Phi}(\delta)(x) \in A^{(m)}$.
Thus $\mu\mi_2(B\mi)=0\,.$ \qed

\noindent
The following proposition is the straighforward corollary of the above lemmas.
\begin{proposition}\label{p56}
$\lim_{m\to\infty} c_{\mu\mi_2}(E)=\frac{3}{2}\,.$
\end{proposition}

\section{The main theorem}
Finally we state and prove our main theorem.
\begin{theorem}\label{tmt}
There exists a minimal topological action $\widetilde{\Phi}$ of the
group $\bZ_2 * \bZ_2 * \bZ_2 * \bZ_2 $ on a compact metric space $Z$
for which there exist ergodic invariant measures $\mu_1$ and $\mu_2$
such that the measurable equivalence relations given by the systems
$(\widetilde{\Phi},Z,\mu_1)$ and $(\widetilde{\Phi},Z,\mu_2)$ are
not orbit equivalent.
\end{theorem}
\proof
By Propositions \ref{p51} and \ref{p56}, for large enough $m$,
the costs of the equivalence relations given by 
$(\widetilde{\Phi}\mi, X_{G\mi} ,\mu\mi_1)$ and 
$(\widetilde{\Phi}\mi ,X_{G\mi} ,\mu\mi_2)$ are
different. By the ergodic decomposition theorem (\cite{KM}, Corollary 18.6),
if $E$ is a Borel-equivalence relation and $\mu$ is an $E$-invariant
measure then
$$c_{\mu}(E)=\int_{\varepsilon I_E} c_e (E) d\nu_\mu (e)\,,$$
where $\varepsilon I_E$ is the space of ergodic
$E$-invariant measures and $\mu=\int e \d\nu_\mu(e)$, where
$\nu_\mu$ is a probability measure.
Hence if $c_{\mu'_1}(E)\neq c_{\mu'_2}(E)$ for some invariant
measures $\mu'_1$ and $\mu'_2$, then there exist ergodic invariant
measures $\mu_1$ and $\mu_2$ for which $c_{\mu_1}(E)\neq c_{\mu_2}(E)$.
Therefore our theorem follows. \qed

\end{document}